\documentstyle[12pt]{article}
\input epsf

\renewcommand{\theequation}{\arabic{section}.\arabic{equation}}
\makeatletter
\@addtoreset{equation}{section}
\makeatother

\newcommand{\beq}{\begin{equation}}
\newcommand{\eeq}{\end{equation}}
\newcommand{\beqa}{\begin{eqnarray}}
\newcommand{\eeqa}{\end{eqnarray}}

\newcommand{\pa}{\partial}

\font\fakeb=msbm10 scaled\magstep1
\def\IC{\mbox{\fakeb C}}

\font\hana=eufm10 scaled\magstep1
\font\shana=eufm10 
\def\HH{\mbox{{\hana H}}}
\def\HH{\mbox{{\hana H}}}
\def\Sn{\mbox{{\hana S}}_N}
\def\sSn{\mbox{{\shana S}}_N}

\newcommand{\ii}{\mbox{i}}

\newcommand{\dd}{\mbox{d}}
\newcommand{\Ker}{\mbox{\rm Ker }}

\newcommand{\nn}{\mbox{\scriptsize\bf n}}

\newtheorem{lem}{Lemma}[section]
\newtheorem{df}[lem]{Definition}
\newtheorem{thm}[lem]{Theorem}

\newtheorem{prop}[lem]{Proposition}

\newcommand{\cx}{\IC\,[x]}
\newcommand{\czz}{\IC\,[z^2]}

\newcommand{\la}{\langle}
\newcommand{\raJ}[1]{\rangle_{\mbox{\rm\scriptsize J}}^{(\beta#1)}}
\newcommand{\raH}[1]{\rangle_{\mbox{\rm\scriptsize H}}^{(\beta#1)}}
\newcommand{\raL}[1]{\rangle_{\mbox{\rm\scriptsize L}}^{(\beta#1)}}
\newcommand{\laa}{\langle\langle}
\newcommand{\raaA}{\rangle\rangle_A}

\newcommand{\Hs}{{\cal H}_{\mbox{s}}}
\newcommand{\tilHs}{\widetilde{\cal H}_{\mbox{s}}}
\newcommand{\Ha}{{\cal H}_A}
\newcommand{\Hb}{{\cal H}_B}
\newcommand{\tilHa}{\widetilde{\cal H}_A}
\newcommand{\tilHb}{\widetilde{\cal H}_B}

\newcommand{\GSs}{\phi_{\mbox{s}}^{(\beta)}}
\newcommand{\GSa}{\phi_{A}^{(\beta)}}
\newcommand{\GSb}{\phi_{B}^{(\beta)}}
\newcommand{\tilphi}{\widetilde{\phi}}

\newcommand{\su}{\mbox{\rm s}}

\newcommand{\nsJ}{E}

\newcommand{\DA}{D^A}
\newcommand{\DB}{D^B}
\newcommand{\hDA}{\hat{D}^A}
\newcommand{\hDB}{\hat{D}^B}
\newcommand{\tilA}{\widetilde{a}}
\newcommand{\tilAdag}{\widetilde{a}^{\dag}}
\newcommand{\tilB}{\widetilde{b}}
\newcommand{\tilBdag}{\widetilde{b}^{\dag}}
\newcommand{\tilhA}{\widetilde{h}^A}
\newcommand{\tilhB}{\widetilde{h}^B}

\newcommand{\hDelJ}{\hat{\Delta}_{\mbox{\rm\scriptsize J}}}
\newcommand{\hDelH}{\hat{\Delta}_{\mbox{\rm\scriptsize H}}}

\newcommand{\sx}{{\cal X}}
\newcommand{\sy}{{\cal Y}}
\newcommand{\hsy}{\hat{\cal Y}}
\newcommand{\syJ}{{\cal Y}_{\mbox{\rm\scriptsize J}}}
\newcommand{\hsyJ}{\hat{\cal Y}_{\mbox{\rm\scriptsize J}}}
\newcommand{\GJ}{G_{\mbox{\rm\scriptsize J}}}
\newcommand{\hGJ}{\hat{G}_{\mbox{\rm\scriptsize J}}}
\newcommand{\tilc}{\widetilde{c}}
\newcommand{\syH}{{\cal Y}_{\mbox{\rm\scriptsize H}}}
\newcommand{\hsyH}{\hat{\cal Y}_{\mbox{\rm\scriptsize H}}}
\newcommand{\GH}{G_{\mbox{\rm\scriptsize H}}}
\newcommand{\hGH}{\hat{G}_{\mbox{\rm\scriptsize H}}}
\newcommand{\sxL}{{\cal X}_{\mbox{\rm\scriptsize L}}}
\newcommand{\syL}{{\cal Y}_{\mbox{\rm\scriptsize L}}}
\newcommand{\hsyL}{\hat{\cal Y}_{\mbox{\rm\scriptsize L}}}
\newcommand{\GL}{G_{\mbox{\rm\scriptsize L}}}
\newcommand{\hGL}{\hat{G}_{\mbox{\rm\scriptsize L}}}

\newcommand{\BmJ}{B^{\mbox{\rm\scriptsize J}}_m}
\newcommand{\BmH}{B^{\mbox{\rm\scriptsize H}}_m}
\newcommand{\BH}[1]{B^{\mbox{\rm\scriptsize H}}_{#1}}
\newcommand{\BmL}{B^{\mbox{\rm\scriptsize L}}_m}
\newcommand{\BL}[1]{B^{\mbox{\rm\scriptsize L}}_{#1}}

\newcommand{\cP}{{\cal P}}
\newcommand{\cPb}{{\cal P}_-^{(\beta)}}
\newcommand{\hw}{\hat{w}}
\newcommand{\hs}{\hat{s}}

\newcommand{\ddaH}{\HH'}
\newcommand{\daH}{\HH_0'}
\newcommand{\ddaHsub}{\widetilde{\HH}'}

\newcommand{\rhoA}{\rho^A}
\newcommand{\rhoB}{\rho^B}
\newcommand{\sigA}{\sigma^A}
\newcommand{\sigB}{\sigma^B}

\newcommand{\Dom}{<_{\mbox{\tiny D}}}
\newcommand{\Bru}{<_{\mbox{\tiny B}}}

\newcommand{\Res}{\mbox{\rm Res}}

\title{Intertwining Operators for a Degenerate\\
 Double Affine Hecke Algebra and\\
 Multivariable Orthogonal Polynomials}
\author{Saburo Kakei\thanks{E-mail: kakei@poisson.ms.u-tokyo.ac.jp}\\[5mm]
{\it Department of Mathematical Sciences, University of Tokyo,}\\
{\it Komaba 3-8-1, Meguro-ku, Tokyo 153, Japan}}
\date{}

\begin{document}
\maketitle

\begin{abstract}
Operators that intertwine representations of a degenerate
version of the double affine Hecke algebra are introduced.
Each of the representations is related to multi-variable orthogonal
polynomials associated with Calogero-Sutherland type models.
As applications, raising operators and shift operators for 
such polynomials are constructed.
\end{abstract}
\bigskip
\bigskip

\section{Introduction}
There are intimate relations between quantum mechanics and 
special functions. 
Wavefunctions for some systems can explicitly be written in terms 
of suitable special functions.
Recent studies on integrable quantum many-particle systems reveal that 
wavefunctions of some special cases can be written in terms of 
multivariable analogue of classical orthogonal polynomials
\cite{BF1,BF2,BF3,vD,K1,K2,So,UW}.
In the present paper, we shall consider orthogonal polynomials associated
with the quantum Calogero models confined
in harmonic potential \cite{Ca1,Ca2,Su1,Y}:
\beqa
\Ha  &=& \frac{1}{2}\sum_{j=1}^{N}
      \left( -\frac{\pa^2}{\pa x_j^2} + x_j^2 \right)
  + \sum_{j<k}\frac{\beta(\beta-1)}{(x_j-x_k)^2},
\label{Ham:CalA}\\
\Hb &=& \frac{1}{2}\sum_{j=1}^{N}
      \left\{ -\frac{\pa^2}{\pa z_j^2} + z_j^2 
             +\frac{\gamma(\gamma-1)}{z_j^2} \right\}
 + \sum_{j<k}\left\{
        \frac{\beta(\beta-1)}{(z_j-z_k)^2} 
        + \frac{\beta(\beta-1)}{(z_j+z_k)^2}\right\},
\label{Ham:CalB} 
\eeqa
where we assume that $\beta$ is a non-negative integer.
The subscripts {\it ``A'', ``B''} signify that
this Hamiltonian is invariant under the action of the Weyl group
of $A_{N-1}$-type or $B_N$-type respectively.
(For the latter convenience, we use the letter ``$z$'' as the 
coordinates of the $B_N$-type model.)
We remark that the model associated with the $C_N$-type Weyl group 
is equivalent to the $B_N$ case, and $D_N$-type model is obtained 
by setting $\gamma=0$.
In these cases, polynomial part of wavefunctions can be regarded as 
multivariable generalization of the Hermite ($A_{N-1}$ case) and Laguerre
($B_N$ case) polynomials and has been studied by several 
authors \cite{BF1,BF2,BF3,vD,K1,K2,So,UW}.

Since special functions are related to representation theory, 
it may be challenging to investigate algebraic aspect of the 
multivariable orthogonal polynomials.
In case of the Macdonald polynomials, it has been revealed that the
algebraic structure behind the polynomials is affine Hecke 
algebras \cite{Ch2,KJr,KN,M2}.
Since the Jack polynomials can be regarded as a degenerate case of 
the Macdonald polynomials, their algebraic structure is a degenerate 
version of affine Hecke algebra \cite{Ch1,Ch3}.

In this paper, we will introduce intertwining operators between 
representations of the degenerate affine Hecke algebra.
From this viewpoint, each of the models (\ref{Ham:CalA}), 
(\ref{Ham:CalB}) corresponds to individual representation of 
the degenerate affine Hecke algebra. 
So far there have been some works concerning representations of 
the algebra associated with (rational) Calogero-type models 
\cite{BF2,BF3,UW}.
However intertwiners between representation spaces have not been 
considered explicitly, though they are important to understand common
algebraic structure of the models.
Using the intertwiners, several results on the Jack polynomials can be
mapped directly to those of the multivariable Hermite and Laguerre 
polynomials.
As applications, we will construct raising operators and shift 
operators for such polynomials.

\section{Dunkl-type operators and multivariable orthogonal polynomials}
\subsection{Jack polynomials and the Sutherland model}
In this subsection, we define our notation and review the theory 
of symmetric and non-symmetric Jack polynomials \cite{M1,O2,KS}.
There are several ways to characterize the Jack polynomials; 
Here we define them as eigenfunctions of some operators.
We note that we restrict ourselves to the case associated with the 
$A_{N-1}$-type Weyl group since we only use such case.

In the paper \cite{D1}, Dunkl has introduced differential-exchange
operators, now called ``Dunkl operators'', which are associated with
root systems. 
For the $A_{N-1}$-type root system, the operators are defined as
\[
\DA_j = \frac{\pa}{\pa x_j}
      + \beta \sum_{k(\neq j)}\frac{1-s_{jk}}{x_j - x_k}
\qquad (j=1,\ldots,N),
\]
where $s_{ij}$ are elements of the symmetric group $\Sn$.
An element $s_{ij}$ acts on functions of $x_1$, $\ldots$, $x_N$ 
as an operator which permutes arguments $x_i$ and $x_j$.
We remark that the operators $D_j$ preserve the space of polynomials
of variables $x_1,\ldots,x_N$ which we denote $\cx$.
These operators satisfy the following properties:
\begin{eqnarray*}
& & [ \DA_i, \DA_j] \;=\; 0,\nonumber\\
& & s_{ij}\DA_j \;=\; \DA_i s_{ij}, \qquad
s_{ij} \DA_k \;=\; \DA_k s_{ij} \quad ( k\neq i,j ),\\[0em]
& & [ \DA_i, x_j ] \;=\; \delta_{ij} 
\left(1+\beta{\displaystyle \sum_{k(\neq i)}}s_{ik}\right)
             -(1-\delta_{ij})\beta s_{ij} .\nonumber
\end{eqnarray*}

Heckman introduced ``global'' Dunkl operators \cite{H2}, 
which are written as $x_j \DA_j$ in our notation.
Heckman's operators do not commute each other.
Cherednik introduced another version of Dunkl operators that
mutually commute \cite{Ch1} (see also \cite{BGHP,KS}):
\beqa
\label{op:ChereA}
\hDA_j & = & x_j \DA_j + \beta \sum_{k(<j)} s_{jk}\nonumber\\
 & = & x_j\frac{\pa}{\pa x_j}
       + \beta \sum_{k(<j)}\frac{x_k}{x_j - x_k} (1-s_{jk})
       + \beta \sum_{k(>j)}\frac{x_j}{x_j - x_k} (1-s_{jk})
       + \beta (j-1) .
\eeqa
The algebra generated by the elements $x_j^{\pm 1}$, $\hDA_j$ and 
$s_{ij}$ is isomorphic to the degenerate double affine Hecke algebra
$\ddaH$ associated with the $A_{N-1}$-type root system \cite{Ch1,Ch3}. 
We remark that the elements $x_j^{\pm 1}$, $\DA_j$ and $s_{ij}$ also
generate $\ddaH$ since $\DA_j$ and $\hDA_j$ are related through
(\ref{op:ChereA}).

We denote by $\daH$ subalgebra of $\ddaH$ generated by $\hDA_j$ and 
$s_{ij}$, which is isomorphic to the degenerate affine Hecke algebra.
We further denote by $\ddaHsub$ subalgebra of $\ddaH$ generated by 
$x_j$, $\hDA_j$ and $s_{ij}$.
In terms of generators, the defining relations are
\begin{eqnarray*}
&& [ \hDA_i , \hDA_j ] = [x_i, x_j] = 0, \nonumber\\
&& s_j ^2 = 1,\qquad s_j s_{j+1} s_j = s_{j+1} s_j s_{j+1},\nonumber\\
&& [s_i, s_j]=0 \quad (|i-j| \neq 1), \nonumber\\
&& x_i s_{ij} = s_{ij} x_j, \qquad 
   x_i s_{jk} = s_{jk} x_i \quad (i \neq j,k),\nonumber\\
&& \hDA_{j+1}s_j - s_j\hDA_j \;=\; \beta,\qquad
   s_j\hDA_{j+1} - \hDA_j s_j \;=\; \beta, \label{rel:daHalg}\\
&& [ s_i, \hDA_j ]=0 \quad (j\neq i,i+1), \nonumber\\
&& [\hDA_i, x_j] = \left\{\begin{array}{lc}
-\beta x_j s_{ij} & (i>j),\\
\displaystyle{
x_i +\beta \left( 
\sum_{k(<i)} x_k s_{ik}+ \sum_{k(>i)} x_i s_{ik} \right) }
 & (i=j), \\
-\beta x_i s_{ij} & (i<j),
\end{array}\right. \nonumber
\end{eqnarray*}
where $s_j=s_{j,j+1}$ ($j=1,\cdots,n-1$) are the simple transpositions.

Since the operators $\hDA_j$ commute each other, they can be
diagonalized simultaneously by suitable choice of bases of $\cx$
\cite{BGHP,O2,KS}.
Such basis is called {\it non-symmetric Jack polynomials}.
To define the non-symmetric Jack polynomials, we first introduce
the ordering $\prec$.

For two pairs $(\lambda, w)$, $(\mu,w')$ where $\lambda$, $\mu$ are 
partitions and $w, w'\in \Sn$,
we define the ordering $\prec$ as follows:
\[
(\mu,w')\prec (\lambda,w) \quad \Longleftrightarrow \quad
\left\{\begin{array}{ll}
(\mbox{\rm i}) & \mu \Dom \lambda,\\
(\mbox{\rm ii}) & \mbox{if } \mu = \lambda \mbox{ then }w'\Bru w,
\end{array}\right.
\]
where $\Dom$ is the dominance ordering for partitions \cite{M1},
and $\Bru$ is the Bruhat ordering for the elements of $\Sn$
(see, for example, \cite{Hu}).
\begin{df}[\cite{BGHP,O2,KS}]
An non-symmetric Jack polynomial $\nsJ^{\lambda}_{w}(x)$,
labeled with the partition $\lambda=(\lambda_1,\ldots,\lambda_N)$ 
and the element $w\in \Sn$, 
is characterized by the following properties:
\begin{enumerate}
\item $\displaystyle \nsJ^{\lambda}_{w}(x) = 
      x^{\lambda}_{w} 
      + \sum_{(\mu,w') \prec (\lambda,w)}
        u^{\lambda \mu}_{w w'} x^{\mu}_{w'}$ ,
\item $\nsJ^{\lambda}_{w}(x)$ is joint eigenfunction for 
      the operators $\hDA_j$,
\end{enumerate}
where we have used the notation 
$x^{\lambda}_{w} = 
x^{\lambda_1}_{w(1)} \cdots x^{\lambda_N}_{w(N)}$.
\end{df}
We note that our definition of the non-symmetric Jack polynomials is 
slightly different from the one in the references cited above. 

Since the action of $\hDA_j$ on monomials $x^{\lambda}_{w}$ are
given by
\beq
\label{eq:ActMon}
\hDA_j x^{\lambda}_{w} = 
(w(\lambda + \beta\delta))_j x^{\lambda}_{w} 
+ \sum_{(\mu,w') \prec (\lambda,w)}
u^{\lambda \mu}_{w w'} x^{\mu}_{w'},
\eeq
with $\delta = (N-1, N-2, \cdots,0)$, 
the action of $\hDA_j$ on the non-symmetric Jack polynomials
can be evaluated as follows \cite{BGHP,O2,KS}:
\beq
\label{eq:eigenV}
\hDA_j \nsJ^{\lambda}_w(x) = 
(w(\lambda + \beta\delta))_j \nsJ^{\lambda}_w(x).
\eeq

From physical viewpoint,
the operators $\hDA_j$ are related to the Sutherland model \cite{Su2}:
\beq
\Hs = -\sum_{j=1}^{N} \frac{\pa^2}{\pa \theta_j^2}
      + \frac{1}{2}\sum_{j<k}
        \frac{\beta(\beta-1)}{
        \sin^2\left[(\theta_j-\theta_k)/2\right]}.
\label{Ham:Suth}
\eeq
To see the relation to the Cherednik operators, we introduce 
``gauge-transformed'' Hamiltonian $\tilHs$:
\begin{eqnarray*}
\lefteqn{\tilHs = \Res\left(\sum_{j=1}^N \left\{ 
\hDA_j - \frac{\beta}{2}(N-1) \right\}^2\right)}
\nonumber\\
&&= \sum_{j=1}^N \left( x_j\frac{\pa}{\pa x_j} \right)^2
       + \beta \sum_{j<k}\frac{x_j+x_k}{x_j-x_k}
       \left( x_j\frac{\pa}{\pa x_j} - x_k\frac{\pa}{\pa x_k}\right)
       + \frac{\beta^2}{12}N(N^2-1),
       \label{eHam:Suth}
\end{eqnarray*}
where $\Res\:X$ means that action of $X$ is restricted to symmetric 
functions of the variables $x_1, \ldots, x_N$.
If we make a kind of gauge transformation and a change of 
variables $x_j = \exp(\ii\theta_j)$, $\tilHs$ reduces to the 
Sutherland Hamiltonian (\ref{Ham:Suth}):
\begin{eqnarray*}
\GSs \circ \tilHs \circ (\GSs)^{-1}
&=& \sum_{j=1}^{N} \left(x_j\frac{\pa}{\pa x_j}\right)^2
      - \beta(\beta-1) \sum_{j<k}
        \frac{2 x_j x_k}{(x_j - x_k)^2}\\
&=& -\sum_{j=1}^{N} \frac{\pa^2}{\pa \theta_j^2}
      + \frac{1}{2}\sum_{j<k}
        \frac{\beta(\beta-1)}{
        \sin^2\left[(\theta_j-\theta_k)/2\right]} = \Hs,
\end{eqnarray*}
where $\GSs(x) = 
       \prod_{j<k}|x_j - x_k|^{\beta}\prod_{j=1}^N x_j^{-\beta(N-1)/2}$
is the ground state wavefunction of the model.
The symmetric Jack polynomials appear as polynomial part of wavefunctions
for excited states.
\begin{df}[\cite{M1}]
\label{def:Jack}
The symmetric Jack polynomials $J_{\lambda}^{(\beta)}(x)$ are 
characterized by the following properties:
\begin{enumerate}
\item $\displaystyle 
J_{\lambda}(x) =  m_{\lambda}(x) + 
\sum_{\mu(\Dom\lambda)}u_{\lambda\mu}m_{\mu}(x)$ ,
\item $J_{\lambda}(x)$ are eigenfunctions of the transformed 
Hamiltonian $\tilHs$,
\end{enumerate}
where $m_{\lambda}$ are the monomial symmetric functions.
\end{df}
The symmetric Jack polynomials are obtained by symmetrizing
$\nsJ^{\lambda}_w$, i.e.,
\[
J_{\lambda}(x) = \frac{1}{\# \Sn^{\lambda}}
\prod_{v\in \sSn}v(\nsJ^{\lambda}_w),
\]
where $\Sn^{\lambda}$ is a subgroup of $\Sn$ that preserve $\lambda$.
This relation follows form the fact that the right hand side
satisfies both of the defining properties
of the Jack polynomials.

As wavefunctions of the Hamiltonian $\Hs$, the following scalar
product is naturally introduced:
\[
\la f(x),g(x) \raJ{} = \oint\cdots\oint f(x)g(x^{-1}) \GSs(x)\GSs(x^{-1})
\frac{\dd x_1}{2\pi\ii x_1}\cdots\frac{\dd x_N}{2\pi\ii x_N},
\]
where the integration contour is the unit circle in the complex plane.
This scalar product can alternatively be written as
\beq
\la f(x),g(x) \raJ{} = (-1)^{\beta N(N-1)/2}
\left[ f \bar{g} (\GSs)^2 \right]_0,
\label{eq:spJ}
\eeq
where $\left[\;\cdot\;\right]_0$ stands for the constant term
and $\bar{g}=g(x^{-1})$.
By a direct calculation, we see that
the operators $\hDA_j$ are self-adjoint with respect to the 
scalar product (\ref{eq:spJ}).
\begin{prop}[\cite{M1}]
The Jack polynomials $J_{\lambda}(x)$ are pairwise orthogonal 
with respect to the scalar product (\ref{eq:spJ}).
\end{prop}
{\it Proof.} We first introduce generating function of 
symmetric commuting operators \cite{BGHP,K1}:
\[
\hDelJ(u) = \prod_{j=1}^N (u + \hDA_j).
\]
If we expand $\hDelJ(u)$ as polynomial in $u$, the coefficients 
form a set of symmetric commuting operators which contains $\tilHs$.
Using (\ref{eq:eigenV}), we can evaluate the action of $\hDelJ(u)$
on the Jack polynomials:
\beq
\label{eq:DelJ}
\hDelJ(u) J_{\lambda}^{(\beta)}(x) = 
\prod_{j=1}^N
\left\{ u + \lambda_{N-j+1} + \beta(j-1) \right\}
J_{\lambda}^{(\beta)}(x).
\eeq
Since all the eigenvalues of $\hDelJ(u)$ are distinct and 
the operator $\hDelJ(u)$ is self-adjoint, we conclude that
the Jack polynomials $J_{\lambda}(x)$ are pairwise orthogonal
with respect to the scalar product (\ref{eq:spJ}).
\hfill $\Box$\\[2mm]

The property below follows from the fact that the Jack polynomials 
form an orthogonal basis of the space of symmetric polynomials
$\cx^{\sSn}$:
\[
\la J_{\lambda}^{(\beta)}(x),m_{\mu}(x)\raJ{} = 0 
\quad \mbox{ for all }\;\mu \Dom \lambda.
\]
One can use this relation instead of the second property of 
Definition \ref{def:Jack}.

\subsection{Multivariable Hermite polynomials and
            $A_{N-1}$-type Calogero model}
We introduce an analogue of the creation and annihilation 
operators:
\[
a^{\dag}_j = \frac{1}{\sqrt{2}} (-\DA_j + x_j), \qquad
a_j = \frac{1}{\sqrt{2}} (\DA_j + x_j).
\]
The commutation relations of these operators are the same as
those of $x_j$ and $\DA_j$ by construction. 
We then make a gauge transformation on $a^{\dag}_j$ and $a_j$:
\begin{eqnarray*}
\tilAdag_j &=& \tilphi^{-1}_A\circ a^{\dag}_j \circ \tilphi_A\\
 &=& \frac{1}{\sqrt{2}}\left(-\frac{\pa}{\pa x_j} 
+2x_j -\beta\sum_{k(\neq j)}\frac{1-s_{jk}}{x_j-x_k}\right),\\
\tilA_j &=& \tilphi^{-1}_A\circ a_j \circ \tilphi_A\\
 &=& \frac{1}{\sqrt{2}}\left(\frac{\pa}{\pa x_j} 
+\beta\sum_{k(\neq j)}\frac{1-s_{jk}}{x_j-x_k}\right),
\end{eqnarray*}
with $\tilphi_A = \prod_{k=1}^N \exp(-x_k^2/2)$.
Since this transformation leaves the commutation relations unchanged, 
we can introduce the following isomorphism:
\[
\rhoA (x_j) = \tilAdag_j, \qquad \rhoA (D_j) = \tilA_j,
\qquad \rhoA (s_{ij}) = s_{ij} .
\]
It should be remarked that this mapping has already been 
appeared implicitly in \cite{UW}, however, treated only as an 
isomorphism of the algebra. 
To construct eigenstates of $\tilHa$, 
we should introduce intertwiner between two representations 
which will be discussed in the followings.

We can obtain a set of commuting operators by applying 
$\rhoA$ to $\hDA_j$:
\[
\tilhA_j = \rhoA (\hDA_j) = \tilAdag_j \tilA_j+ \beta \sum_{k(<j)}s_{jk}.
\]
The mapping $\rhoA$ gives another 
representation of $\ddaHsub$ on $\cx$.
We then introduce intertwining operator $\sigA$, which is a linear
operator on $\cx$ such that
\[
\sigA\left(f(x_1,\ldots,x_N)\right) = 
f(\tilAdag_1,\ldots,\tilAdag_N)\cdot 1
\quad\mbox{for all}\quad f(x_1,\ldots,x_N)\in\cx.
\]
The intertwiner $\sigA$ enjoys the following property.
\begin{thm}
\label{thm:intertwinerA}
$\sigA (Q f(x))= \rhoA(Q)\sigA(f(x))$ for all $Q\in\daH$, 
$f(x)\in\cx$.
\end{thm}
{\it Proof.} 
Since both $Q$ and $f(x)$ are elements of $\ddaHsub$,
it suffices to prove
$\sigA (P\cdot 1)= \rhoA(P)\cdot 1$ for all $P\in\ddaHsub$.
We then note that every element $P$ of $\ddaHsub$ can be
represented in the following form:
\beq
\label{eq:P}
P = \sum_{\nn}\sum_{w\in\sSn}p_{\nn,w}(x)
(\hDA)^{n_1 }\cdots(\hDA_N)^{n_N} w,
\eeq
where $p_{\nn,w}(x)$ are some polynomials.
Considering the action of (\ref{eq:P}) on $1$, we have
\[
P\cdot 1 = \sum_{\nn (n_1=0)}\sum_{w\in\sSn}p_{\nn,w}(x)
\beta^{n_2}\cdots((N-1)\beta)^{n_N},
\]
since $w\cdot 1=1$ for all $w\in\Sn$ and 
$\hDA_j\cdot 1 =\beta(j-1)$ for all $j$.
On the other hand, applying $\rhoA$ to (\ref{eq:P}), we have
\[
\rhoA(P) = \sum_{\nn}\sum_{w\in\sSn}
p_{\nn,w}(\tilAdag) (\tilhA)^{n_1 }\cdots(\tilhA_N)^{n_N} w.
\]
Since $\tilhA_j\cdot 1=\beta(j-1)$ for all $j$, we conclude that
$\sigA (P\cdot 1)= \rhoA(P)\cdot 1$ for all $P\in\ddaH$.
\hfill $\Box$\\[2mm]

The representation $\rhoA$ is related to the $A_{N-1}$-type 
Calogero model.
If we define $\tilHa$ as 
\begin{eqnarray*}
\tilHa &=& \Res\left(\sum_{j=1}^N \tilhA_j\right) 
- \frac{\beta}{2}N(N-1) \nonumber\\
 &=& \frac{1}{2}\sum_{j=1}^n \left( -\frac{\pa^2}{\pa x_j^2} 
+ 2 x_j \frac{\pa}{\pa x_j} \right) - \beta \sum_{j < k}
\frac{1}{x_j - x_k} \left(\frac{\pa}{\pa x_j}-\frac{\pa}{\pa x_k}
\right),
\end{eqnarray*}
we can obtain the $A_{N-1}$-type Calogero Hamiltonian (\ref{Ham:CalA})
via gauge transformation:
\[
\Ha = \GSa \circ\tilHa\circ (\GSa)^{-1}
+ \frac{N}{2} + \frac{\beta}{2}N(N-1),
\]
with 
$\GSa = \prod_{j<k}|x_j - x_k|^{\beta}\prod_{j=1}^N \exp (-x_j^2/2)$
ground state wavefunction.

We then introduce scalar product for this case:
\beq
\la f,g \raH{} = \int_{-\infty}^{\infty}\cdots\int_{-\infty}^{\infty}
f(x)g(x) (\GSa)^2 \dd x_1\cdots\dd x_N 
\label{eq:spH}
\eeq
By a direct calculation, we see that
the operator $\tilAdag_j$ is adjoint of $\tilA_j$ with respect to 
the scalar product (\ref{eq:spH}) for all $j=1,\ldots,N$.
Note that $x_j$ ($=(\rhoA)^{-1}(\tilAdag_j)$) is not adjoint of 
$D_j$($=(\rhoA)^{-1}(\tilA_j)$) for the Jack case.

Multivariable Hermite polynomials are defined by using this scalar
product \cite{BF1,vD}. In fact, the definition in \cite{BF1} and
that in \cite{vD} are slightly different. 
Here we shall follow \cite{vD}:
\begin{df}[\cite{vD}]
\label{def:multiH}
Multivariable Hermite polynomials $H^{(\beta)}_{\lambda}(x)$ are 
characterized by the following properties:
\begin{enumerate}
\item $\displaystyle 
H^{(\beta)}_{\lambda}(x) =  m_{\lambda}(x) + 
\sum_{\mu(\Dom\lambda)}u_{\lambda\mu}^A m_{\mu}(x)$,
\item $\displaystyle \la H^{(\beta)}_{\lambda}(x),m_{\mu}(x)\raH{} = 0 
\quad \mbox{ for all }\;\mu \Dom \lambda$.
\end{enumerate}
\end{df}

Using the intertwiner $\sigA$, we can construct an operator representation
of $H^{(\beta)}_{\lambda}(x)$.
\begin{prop}[\cite{K1,UW}]
\label{prop:opH}
Multivariable Hermite polynomials $H^{(\beta)}_{\lambda}(x)$ are
related to the Jack polynomials as follows:
\[
H^{(\beta)}_{\lambda}(x) = 
2^{-|\lambda|/2}\sigA(J^{(\beta)}_{\lambda}(x))
= 2^{-|\lambda|/2}J^{(\beta)}_{\lambda}(\tilAdag)\cdot 1.
\]
\end{prop}
{\it Proof.} 
We can easily see that $2^{|\lambda|/2}
J_{\lambda}(\tilAdag_1,\ldots,\tilAdag_N)\cdot 1$
satisfy the condition (i) of Definition \ref{def:multiH}.
Hence it suffices to show (ii).
Applying $\sigA$ to (\ref{eq:DelJ}), we have
\[
\hDelH(u) J^{(\beta)}_{\lambda}(\tilAdag_1,\ldots,\tilAdag_N)\cdot 1 = 
\prod_{j=1}^N
\left\{ u + \lambda_{N-j+1} + \beta(j-1) \right\}
J^{(\beta)}_{\lambda}(\tilAdag_1,\ldots,\tilAdag_N)\cdot 1,
\]
where we denote $\hDelH(u)=\rhoA(\hDelJ(u))=\prod_{j=1}^N (u + \tilhA_j)$.
Since all the eigenvalues of $\hDelH(u)$ are distinct and 
the operator $\hDelH(u)$ are self-adjoint with respect to the scalar 
product (\ref{eq:spH}),
we conclude that the polynomials
$J^{(\beta)}_{\lambda}(\tilAdag)\cdot 1$ 
are orthogonal with respect to the scalar product (\ref{eq:spH}).
On the other hand, one may know that the polynomials
$J^{(\beta)}_{\lambda}(\tilAdag)\cdot 1$ 
form an orthogonal basis of $\cx^{\sSn}$ by considering the leading term. 
It follows that
$\la J^{(\beta)}_{\lambda}(\tilAdag)\cdot 1 ,m_{\mu}\raJ{} = 0$ for all 
$\mu \Dom \lambda$, which proves the theorem.
\hfill $\Box$\\[2mm]

\noindent
It should be noted that Ujino and Wadati \cite{UW} have shown that 
$J^{(\beta)}_{\lambda}(\tilAdag)\cdot 1$ diagonalize the first two
of the family of commuting operators that contains $\tilHa$.
The proof given here is essentially the same as that given in \cite{K1}.

The scalar product $\la\cdot,\cdot\raH{}$ induces another scalar 
product on $\cx$:
\[
\label{thirdSP1}
\laa f(x),g(x) \raaA = \la f(\tilAdag)\cdot 1, g(\tilAdag)\cdot 1 \raH{}.
\]
This gives another example of scalar product
which makes the Jack polynomials orthogonal.
On the other hand, Dunkl \cite{D2} introduced the scalar product
$\left[ f(\hDA)g(x)\right]_0$.
These scalar products coincide up to a constant factor:
\begin{eqnarray*}
\lefteqn{
\laa f(x),g(x) \raaA = \la 1, f(\tilA)g(\tilAdag)\cdot 1 \raH{}}\\
&& = \la 1,1\raH{} \left[ f(\tilA)g(\tilAdag)\cdot 1 \right]_0
= \la 1,1\raH{} \left[ f(\hDA)g(x)\right]_0.
\end{eqnarray*}
We shall evaluate the value $\la 1,1\raH{}$ in section \ref{sec:shiftH}.
(See Proposition \ref{prop:normH} below.)

\subsection{Multivariable Laguerre polynomials and
            $B_N$-type Calogero model}
Dunkl operators associated with the $B_N$-type root system are 
defined as follows \cite{D1,Y}:
\beq
\DB_j = \frac{\pa}{\pa z_j}
+ \beta \sum_{k(\neq j)}\left(\frac{1-s_{jk}}{z_j - z_k} 
     + \frac{1-t_j t_k s_{jk}}{z_j + z_k} \right)
     + \gamma \frac{1-t_j}{z_j},
\label{eq:B-Dunkl}
\eeq
where $s_{jk}$ and $t_j$ are elements of the $B_N$-type Weyl group.
An element $s_{ij}$ acts as same as in the $A_{N-1}$-case and 
$t_j$ acts as sign-change, i.e. replaces the coordinate $z_j$ by $-z_j$.
Commutation relations of the $B_N$-type Dunkl operators are
\[
\begin{array}{l}
\displaystyle
[ \DB_i, \DB_j] = 0,\\
\displaystyle
[ \DB_i, z_j ] = 
      \delta_{ij}\left\{1+ \beta
         \sum_{k(\neq i)}(s_{ik}+t_i t_k s_{ik})
         + 2\gamma t_j \right\}\\
\qquad\qquad\quad -(1-\delta_{ij})\beta (s_{ij}-t_i t_k s_{ik}),\\
s_{ij}\DB_j = \DB_i s_{ij}, \qquad
s_{ij} \DB_k = \DB_k s_{ij} \quad ( k\neq i,j ),\\
t_j\DB_j = -\DB_i t_j, \qquad
t_j \DB_k = \DB_k t_j \quad ( k\neq j ).
\end{array}
\]

We then define Cherednik-type commuting operators 
associated with (\ref{eq:B-Dunkl}):
\[
\hDB_j = z_j \DB_j 
       + \beta \sum_{k(<j)}(s_{jk} + t_j t_k s_{ik}).
\]
Note that the operators $\hDB_j$ do {\it not} coincide with 
the Cherednik operators associated with the $B_N$-type Weyl group.
\begin{lem}
All of the operators $\hDB_j$, $s_{ij}$, $t_j$ and $z_j^2$ preserve 
$\IC[z_1^2,\ldots,z_N^2]$.
\end{lem}
{\it Proof.} 
Only $\hDB_j$ need to prove.
We introduce the notation $\Res^{(t)}(X)$ which means the
action of the operator $X$ is restricted to the functions with 
the symmetry $t_j f(z) = f(z)$.
Under this restriction, the action of the operator $\hDB_j$ is 
reduced to the following form:
\begin{eqnarray}
\Res^{(t)}(\hDB_j) & = & z_j\frac{\pa}{\pa z_j}
  + 2\beta \sum_{k(<j)}\frac{z_k^2}{z_j^2 - z_k^2}(1-s_{jk})
\nonumber\\
& & \qquad + 2\beta \sum_{k(>j)}\frac{z_j^2}{z_j^2 - z_k^2} (1-s_{jk})
  + 2\beta (j-1).
\label{op:ResChereB}
\end{eqnarray}
Comparing (\ref{op:ResChereB}) with (\ref{op:ChereA}), we find that
$\Res^{(t)}(\hDB_j)$ is equivalent to $2\hDA_j$ if we make a change 
of the variables $x_j=z_j^2/2$.
Since $\hDA_j$ preserve $\cx$, the operators $\hDB_j$ preserve $\czz$.
\hfill $\Box$\\[2mm]
From these facts, we can define representation $\iota$ of 
$\ddaHsub$ on $\czz$:
\[
\iota (x_j) = \frac{1}{2}z_j^2, \qquad
\iota (\hDA_j) = \frac{1}{2}\hDB_j, \qquad
\iota (s_{ij}) = s_{ij}.
\]

We now introduce creation and annihilation operators for the $B_N$ case:
\[
b^{\dag}_j = \frac{1}{\sqrt{2}} (-\DB_j + z_j), \qquad
b_j = \frac{1}{\sqrt{2}} (\DB_j + z_j).
\]
The commutation relations of these operators are the same as
those of $z_j$ and $\DB_j$ by construction. 
We then make a gauge transformation on $b^{\dag}_j$ and $b_j$:
\begin{eqnarray*}
\tilBdag_j &=& \tilphi_B^{-1}\circ b^{\dag}_j \circ \tilphi_B\nonumber\\
	 &=& \frac{1}{\sqrt{2}}\left\{-\frac{\pa}{\pa z_j} 
+2z_j -\beta\sum_{k(\neq j)}\left(\frac{1-s_{jk}}{z_j - z_k} 
     + \frac{1-t_j t_k s_{jk}}{z_j + z_k} \right)
     + \gamma \frac{1-t_j}{z_j}\right\},\\
\tilB_j &=& \tilphi_B^{-1}\circ b_j \circ \tilphi_B\nonumber\\
 &=& \frac{1}{\sqrt{2}}\left\{\frac{\pa}{\pa z_j} 
+\beta\sum_{k(\neq j)}\left(\frac{1-s_{jk}}{z_j - z_k} 
     + \frac{1-t_j t_k s_{jk}}{z_j + z_k} \right)
     + \gamma \frac{1-t_j}{z_j}\right\},
\end{eqnarray*}
with $\tilphi_B = \prod_{k=1}^N \exp(-z_k^2/2)$.
Since this transformation leaves the commutation relations unchanged, 
we can define the following algebra isomorphism:
\[
\kappa (x_j) = \tilBdag_j, \qquad \kappa (\DB_j) = \tilB_j,
\qquad \kappa (s_{ij}) = s_{ij}, \qquad \kappa (t_j) = t_j.
\]

We then define the operators $\tilhB_j$ as follows:
\[
\tilhB_j = \kappa(\hDB_j) = 
\tilBdag_j \tilB_j 
       + \beta \sum_{k(<j)}(s_{jk} + t_j t_k s_{ik}).
\]
\begin{lem}
The operators $\tilhB_j$ and $(\tilBdag_j)^2$ preserve 
$\IC[z_1^2,\ldots,z_N^2]$.
\end{lem}
{\it Proof.} 
Since the operators $\DB_j$ preserve $\IC[z_1,\ldots,z_N]$, 
it is clear that both $\tilhB_j$ and $(\tilBdag_j)^2$ also
preserve $\IC[z_1,\ldots,z_N]$. 
Then it suffices to prove 
$[t_i,\tilhB_j]=[t_i,(\tilBdag_j)^2]=0$ for all $i,j$,
which can be proved by a direct calculation.
\hfill $\Box$\\[2mm]

Using both $\iota$ and $\kappa$, we introduce another representation
of $\ddaHsub$ on $\czz$:
\[
\rhoB(x_j) = \kappa (\iota(x_j)) = \frac{1}{2}(\tilBdag_j)^2, \qquad
\rhoB(\hDA_j) = \kappa(\iota(\hDA_j))
= \frac{1}{2}\tilhB_j, \qquad
\rhoB(s_{ij}) = s_{ij}.
\]
We introduce a linear map of $\cx$ to $\czz$ by using $\rhoB$:
\[
\sigB (f(x_1,\ldots,x_N)) = 
f((\tilBdag_1)^2/2,\ldots,(\tilBdag_N)^2/2)\cdot 1
\quad\mbox{for all}\quad f(x_1,\ldots,x_N)\in\cx.
\]
As in the $A_{N-1}$-case, the intertwiner $\sigB$ enjoys the 
following property.
\begin{thm}
$\sigB (Q f(x))= \rhoB(Q)\sigB(f(x))\cdot 1$ 
for all $Q\in\daH$, $f(x_1,\ldots,x_N)\in\cx$.
\end{thm}
Proof is given in the same fashion as Theorem \ref{thm:intertwinerA},
so we omit details.

The operators $\tilBdag_j$ and $\tilB_j$ are related to the $B_N$-type 
Calogero Hamiltonian (\ref{Ham:CalB});
If we define $\tilHb$ as 
\begin{eqnarray*}
\tilHb &=& \Res\left(\sum_{j=1}^N \tilhB_j\right) 
-\beta N(N-1) \nonumber\\
 &=& \frac{1}{2}\sum_{j=1}^n \left( -\frac{\pa^2}{\pa z_j^2} 
+2 z_j \frac{\pa}{\pa z_j}-\frac{2\gamma}{z_j}\frac{\pa}{\pa z_j}
\right) - 2\beta \sum_{j < k}
\frac{1}{z_j^2 - z_k^2} 
\left(z_j\frac{\pa}{\pa z_j}-z_k\frac{\pa}{\pa z_k}\right),
\end{eqnarray*}
we can obtain the Hamiltonian (\ref{Ham:CalB})
via gauge transformation:
\[
\Hb = \GSb \circ\tilHb\circ (\GSb)^{-1}
+ \left(\frac{1}{2}+\gamma\right)N + \beta N(N-1),
\]
with 
$\GSb = \prod_{j<k}|z_j^2 - z_k^2|^{\beta}\prod_{j=1}^N 
|z_j|^{\gamma}\exp (-z_j^2/2)$
ground state wavefunction of (\ref{Ham:CalB}).

Scalar product associated with this model is
\beq
\label{eq:spL}
\la f(z),g(z)\raL{} = \int_{-\infty}^{\infty}\cdots\int_{-\infty}^{\infty}
f(z)g(z)(\GSb)^2\dd z_1\cdots\dd z_N .
\eeq
By a direct calculation, we can show that
the operator $\tilBdag_j$ is adjoint of 
$\tilB_j$ with respect to the scalar product (\ref{eq:spL}),
and hence the operator $\tilhB_j$
is self-adjoint for all $j=1,\ldots,N$.

Now we define multivariable Laguerre polynomials \cite{vD}.
\begin{df}[\cite{vD}]
\label{def:multiL}
Multivariable Laguerre polynomials $L^{(\beta)}_{\lambda}(z)$ are 
characterized by the following properties:
\begin{enumerate}
\item $\displaystyle 
L^{(\beta)}_{\lambda}(z) =  m_{\lambda}(z^2) + 
\sum_{\mu(\Dom\lambda)}u_{\lambda\mu}m_{\mu}(z^2)$,
\item $\displaystyle \la L^{(\beta)}_{\lambda}(z),m_{\mu}(z^2)\raL{} = 0 
\quad \mbox{ for all }\;\mu \Dom \lambda$.
\end{enumerate}
\end{df}
We can construct an operator representation
of $L^{(\beta)}_{\lambda}(z)$ by using the intertwiner $\sigB$.
\begin{prop}[\cite{K2}]
\label{prop:opL}
Multivariable Laguerre polynomials $L^{(\beta)}_{\lambda}(z)$ are
related to the Jack polynomials as follows:
\[
L^{(\beta)}_{\lambda}(z) = \sigB(J^{(\beta)}_{\lambda}(x))
= J^{(\beta)}_{\lambda}
((\tilBdag)^2/2)\cdot 1.
\]
\end{prop}
One can prove this statement in the same way as Proposition 
\ref{prop:opH}, so we omit details.

\section{Construction of raising operators}
As is shown in the last section, the multivariable Hermite and
Laguerre polynomials are expressed in terms of  the Jack polynomials
whose arguments are Dunkl-type operators. 
Some properties of the multivariable Hermite and Laguerre polynomials 
are obtained directly from those of the Jack polynomials simply by 
applying $\rhoA$ or $\rhoB$.
As an example, we will construct raising operators for the polynomials.

Lapointe and Vinet constructed raising operators for the Jack 
polynomials \cite{LV}.
Using their raising operators, they obtained Rodorigues-type formula
for the Jack polynomials.
Raising operators for the multivariable Hermite polynomials have been
constructed by Ujino and Wadati \cite{UW}.
The raising operators above are constructed by the use of 
Dunkl operators of Heckman-type (non-commutative).

On the other hand, Kirillov and Noumi gave another expression of
raising operators by using Cherednik operators \cite{KN}. 
In our notation, their raising operators are expressed as the
following form:
\begin{eqnarray*}
\BmJ &=& 
\sum_{k_1<\cdots<k_m}
x_{k_1}x_{k_2}\cdots x_{k_m}
(\hDA_{k_1}+\beta(2-k_1))\\
&& \qquad\qquad\times(\hDA_{k_2}+\beta(3-k_2))
\cdots(\hDA_{k_m}+\beta(m-k_m+1)).
\end{eqnarray*}
We recall a important property of these operators.
\begin{thm}[\cite{KN}]
\label{thm:RaiseJ}
Action of the operators $\BmJ\in\ddaHsub$ on the Jack polynomials are
given by 
\[
\BmJ J_{\lambda}^{(\beta)}(x) =
\prod_{j=1}^m (\lambda_j + \beta(m-j+1))
J_{\lambda+(1^m)}^{(\beta)}(x),
\]
where $\lambda +(1^m) = (\lambda_1+1,\ldots,\lambda_N+1)$.
\end{thm}

Applying $\sigA$ or $\sigB$ to $\BmJ$, we obtain raising operators
for the Hermite-case or the Laguerre-case respectively:
\begin{eqnarray*}
\BmH &=& \sum_{k_1<\cdots<k_m}
\tilAdag_{k_1}\tilAdag_{k_2}\cdots\tilAdag_{k_m}
(\tilhA_{k_1}+\beta(2-k_1))\\
&& \qquad\qquad\times (\tilhA_{k_2}+\beta(3-k_2))
\cdots(\tilhA_{k_m}+\beta(m-k_m+1)),\\
\BmL &=& \sum_{k_1<\cdots<k_m}
\tilBdag_{k_1}\tilBdag_{k_2}\cdots\tilBdag_{k_m}
(\tilhB_{k_1}+\beta(2-k_1))\\
&& \qquad\qquad\times (\tilhB_{k_2}+\beta(3-k_2))
\cdots(\tilhB_{k_m}+\beta(m-k_m+1)).
\end{eqnarray*}
Form the theorem \ref{thm:RaiseJ} and
the propositions \ref{prop:opH}, \ref{prop:opL}, 
it immediately follows that:
\begin{prop}
\begin{enumerate}
\item
$\BmH H_{\lambda}^{(\beta)}(x) =
2^{-m/2}\prod_{j=1}^m (\lambda_j + \beta(m-j+1))
H_{\lambda+(1^m)}^{(\beta)}(x)$,
\item
$\BmL L_{\lambda}^{(\beta)}(z) =
\prod_{j=1}^m (\lambda_j + \beta(m-j+1))
L_{\lambda+(1^m)}^{(\beta)}(z)$.
\end{enumerate}
\end{prop}

Applying the raising operators repeatedly, one can obtain 
Rodorigues-type formulas for the multivariable Hermite and Laguerre
polynomials:
\begin{eqnarray*}
H_{\lambda}^{(\beta)}(x) &=& 2^{|\lambda|/2}
\prod_{(i,j)\in\lambda}(\lambda_i-j + \beta(\lambda'_j-i+1))^{-1}
(\BH{N})^{\lambda_N} (\BH{N-1})^{\lambda_{N-1}-\lambda_N}
\cdots (\BH{1})^{\lambda_1-\lambda_2}\cdot 1,\\
L_{\lambda}^{(\beta)}(z) &=& 
\prod_{(i,j)\in\lambda}(\lambda_i-j + \beta(\lambda'_j-i+1))^{-1}
(\BL{N})^{\lambda_N} (\BL{N-1})^{\lambda_{N-1}-\lambda_N}
\cdots (\BL{1})^{\lambda_1-\lambda_2}\cdot 1,
\end{eqnarray*}
where $\lambda'=(\lambda'_1,\lambda'_2,\ldots)$ is the conjugate 
partition to $\lambda$.

\section{Construction of shift operators}
In this section, we construct shift operators for the multivariable 
Hermite and Laguerre polynomials. 
Each of such shift operators are related to one of the scalar 
products (\ref{eq:spJ}), (\ref{eq:spH}), (\ref{eq:spL}) 
(see the theorems \ref{thm:dualityJ}, \ref{thm:dualityH} and 
\ref{thm:dualityL} below).
However properties related to scalar product cannot be obtained 
simply by applying $\rhoA$ or $\rhoB$.
It require a little more effort to construct shift operators.

\subsection{Shift operators for the Jack polynomials}
In this subsection, we review the method of constructing shift operators
for the Jack polynomials by the use of the Cherednik operators.
Our method is based on the lecture note of Kirillov Jr. \cite{KJr};
All results given in this section can be obtained by taking limiting
procedure on those of \cite{KJr}.
However, all proofs given here are algebraic and we avoid using
limiting procedure so that we can apply the results to the 
Hermite and Laguerre cases.

Consider elements of $\daH$ that have the following forms:
\begin{eqnarray*}
\sx &=& \prod_{i<j}(x_i-x_j), \nonumber\\
\syJ &=& \prod_{i<j}(\beta -\hDA_i + \hDA_j),\label{eq:sxsy}\\
\hsyJ &=& \prod_{i<j}(-\beta -\hDA_i + \hDA_j).\nonumber
\end{eqnarray*}
We note that these operators preserve $\cx$.
From the defining relations of $\daH$, we know that
\beqa
(s_j+1)(-\beta -\hDA_j+\hDA_{j+1}) &=& (-\beta -\hDA_{j+1}+\hDA_j)(s_j-1),
\nonumber\\
(s_j-1)(\beta -\hDA_j+\hDA_{j+1}) &=& (\beta -\hDA_{j+1}+\hDA_j)(s_j+1),
\label{eq:sDD}\\
s_j (c -\hDA_j+\hDA_k)(c -\hDA_{j+1}+\hDA_k)
&=& (c -\hDA_j+\hDA_k)(c -\hDA_{j+1}+\hDA_k) s_j, \nonumber
\eeqa
with $c$ arbitrary constant and $k\neq j,j+1$.
Then, if we define $\cx^{\sSn}$ and $\cx^{-\sSn}$ as
\begin{eqnarray*}
\cx^{\sSn} &=& \{ f(x)\in\cx \;|\; (s_j -1)f(x)=0 \}, \\
\cx^{-\sSn} &=& \{ f(x)\in\cx \;|\; (s_j +1)f(x)=0 \},
\end{eqnarray*}
we see that
\beq
\label{eq:xyy}
\sx \in \cx^{-\sSn}, \qquad
\syJ\left(\cx^{\sSn}\right) = \cx^{-\sSn}, \qquad
\hsyJ\left(\cx^{-\sSn}\right) = \cx^{\sSn}.
\eeq

We now introduce shift operators for the Jack polynomials:
\[
\GJ = \sx^{-1}\syJ, \qquad \hGJ = \hsyJ\sx.
\]
The operators $\GJ$ and $\hGJ$ enjoy the following properties:
\begin{lem}[\cite{KJr}]
\label{lem:actGJ}
\begin{enumerate}
\item
$\GJ J_{\lambda+\delta}^{(\beta)}$, $\hGJ J_{\lambda}^{(\beta+1)}$ 
$\in\cx^{\sSn}$.
\item
$\GJ J_{\lambda+\delta}^{(\beta)} = c_{\lambda}^{(\beta+1)} m_{\lambda} + $
``lower terms'' with respect to $\Dom$,\\
with $c_{\lambda}^{(\beta+1)} = \prod_{i<j}
\left\{\lambda_{N-j+1}-\lambda_{N-i+1}+j-i+\beta(j-i-1)\right\}$.
\item
$\hGJ J_{\lambda}^{(\beta+1)} = \tilc_{\lambda}^{(\beta+1)} 
 m_{\lambda+\delta} + $ ``lower terms'' with respect to $\Dom$,\\
with $\tilc_{\lambda}^{(\beta+1)} = \prod_{i<j}
\left\{\lambda_{N-j+1}-\lambda_{N-i+1}+j-i+\beta(j-i+1)\right\}$.
\end{enumerate}
\end{lem}
{\it Proof.} 
(i) Follows from (\ref{eq:xyy}).

(ii) 
For the longest element $w_0$ of $\Sn$, i.e. 
$w_0(j)=N-j+1$, equation (\ref{eq:ActMon}) reduces to
\[
\hDA_j x^{\lambda}_{w_0} = 
(\lambda_{N-j+1} + \beta(j-1)) x^{\lambda}_{w_0} 
+ \sum_{(\mu,w') \prec (\lambda,w_0)}
u^{\lambda \mu}_{w_0 w'} x^{\mu}_{w'}.
\]
Using this relation, we can calculate the action of $\syJ$:
\begin{eqnarray}
\syJ J_{\lambda+\delta}^{(\beta)} &=& \prod_{i<j}(\beta -\hDA_i + \hDA_j)
(x^{\lambda+\delta}_{w_0}+\mbox{ ``lower terms'' with respect to $\prec$})
\nonumber\\
&=& c_{\lambda}^{(\beta+1)}x^{\lambda+\delta}_{w_0} 
    +\mbox{ ``lower terms'' with respect to $\prec$}.
\label{eq:actJ}
\end{eqnarray}
On the other hand, (\ref{eq:xyy}) implies that 
$\syJ J_{\lambda+\delta}^{(\beta)}$ is divisible by $\sx$.
Together with (i), this concludes the proof.

(iii) Can also be proved in similar way.
\hfill $\Box$\\[2mm]

The following theorem implies that $\hGJ$ is, in a sense, 
adjoint of $\GJ$:
\begin{thm}[\cite{KJr}]
\label{thm:dualityJ}
For $f, g \in \cx^{\sSn}$, 
$\displaystyle
\la \GJ f, g\raJ{+1} = \la f, \hGJ g\raJ{}$.
\end{thm}
To prove this theorem, we introduce symmetrizer $\cP_+$ and 
anti-symmetrizer $\cP_-$ as 
\[
\cP_+ = \frac{1}{\# \Sn}\sum_{w\in \sSn}w,\qquad
\cP_- = \frac{1}{\# \Sn}\sum_{w\in \sSn}(-1)^{l(w)}w,
\]
where $l(w)$ is the length of the element $w$.
We further prepare a lemma.
\begin{lem}[\cite{KJr}]
\label{lem:PY2}
$\displaystyle
\cP_- (\syJ - \hsyJ) = \sum_j \hat{g}_j(\hDA_1,\ldots,\hDA_N)(s_j-1)$
for some $\hat{g}_j(x_1,\cdots,x_N)\in\cx$.
\end{lem}
It should be remarked that this lemma is a degenerate version of 
a lemma given in \cite{KJr}.
We will give a proof in Appendix A for reader's convenience.

Now we go back to the proof of Theorem \ref{thm:dualityJ}.\\[2mm]
{\it Proof of Theorem \ref{thm:dualityJ}.} 
It is clear that $\cP_+$ does not affect 
constant term of polynomials.
Then, for $f,g \in \cx^{\sSn}$, we know that
\begin{eqnarray*}
\la \GJ f, g \raJ{+1} &=& (-1)^{(\beta+1)N(N-1)/2}
\left[ (\sx^{-1}\syJ f)\bar{g}(\phi_{\su}^{(\beta+1)})^2 \right]_0\\
&=& (-1)^{(\beta+1)N(N-1)/2}\left[ \cP_+ \left((\syJ f)\bar{g}\sx 
(\GSs)^2 \right)\right]_0\\
&=& (-1)^{\beta N(N-1)/2}\left[ \cP_- (\syJ f)\bar{g}\bar{\sx}
    (\GSs)^2 \right]_0.
\end{eqnarray*}
From Lemma \ref{lem:PY2}, we see that
$\cP_- (\syJ - \hsyJ) f =0$ for all $f\in\cx^{\sSn}$.
Hence we can replace $\syJ$ by $\hsyJ$:
\begin{eqnarray*}
\la \GJ f, g \raJ{+1} &=& (-1)^{\beta N(N-1)/2}
\left[ \cP_- (\hsyJ f)\bar{g}\bar{\sx}(\GSs)^2 \right]_0\\
&=& \la \hsyJ f, \sx g \raJ{} = \la f, \hGJ g \raJ{}.
\end{eqnarray*}
In the last equality, we have used the self-adjointness of $\hDA_j$.
\hfill $\Box$\\[2mm]

Using Theorem \ref{thm:dualityJ}, we can evaluate the action of 
$\GJ$ and $\hGJ$ on the Jack polynomials.
\begin{prop}[\cite{KJr}]
\label{prop:SrelJ}
$\GJ J_{\lambda+\delta}^{(\beta)} = 
c_{\lambda}^{(\beta+1)} J_{\lambda}^{(\beta+1)}$,
$\hGJ J_{\lambda}^{(\beta+1)} = \tilc_{\lambda}^{(\beta+1)} 
J_{\lambda+\delta}^{(\beta)}$,
where the constants $c_{\lambda}^{(\beta)}$ and 
$\tilc_{\lambda}^{(\beta)}$ are defined in Lemma \ref{lem:actGJ}.
\end{prop}
{\it Proof.} 
Assume $\mu\Dom\lambda$. Then we have
\begin{eqnarray*}
\lefteqn{
\la \GJ J_{\lambda+\delta}^{(\beta)}, m_{\mu}\raJ{}
= \la J_{\lambda+\delta}^{(\beta)}, \hGJ m_{\mu}\raJ{}}\\
&& = \la J_{\lambda+\delta}^{(\beta)}, m_{\mu+\delta}+
\mbox{ lower terms }\raJ{} = 0,
\end{eqnarray*}
where we have used Lemma \ref{lem:actGJ} and Theorem \ref{thm:dualityJ}.
From this fact, along with Lemma \ref{lem:actGJ} (ii), we see that 
$\GJ J_{\lambda+\delta}^{(\beta)}$ coincides with
$c_{\lambda}^{(\beta+1)}J_{\lambda}^{(\beta+1)}$.
The latter can be proved in similar way.
\hfill $\Box$\\[2mm]

With these preliminaries, it is possible to derive the following
result.
\begin{prop}[\cite{M1}]
\beq
\label{eq:NormJack}
\la J_{\lambda}^{(\beta)},J_{\lambda}^{(\beta)}\raJ{}
= N! \prod_{k=1}^{\beta}\prod_{i<j}
\frac{\lambda_i - \lambda_j -k + \beta(j-i+1)}{
\lambda_i - \lambda_j +k + \beta(j-i-1)}.
\eeq
\end{prop}
{\it Proof.} 
From Proposition \ref{prop:SrelJ}, it follows that
\begin{eqnarray*}
\lefteqn{\la J_{\lambda}^{(\beta+1)},J_{\lambda}^{(\beta+1)}\raJ{+1}
= \frac{1}{c_{\lambda}^{(\beta+1)}}
\la G J_{\lambda+\delta}^{(\beta)},J_{\lambda}^{(\beta+1)}\raJ{+1}}\\
&& = \frac{1}{c_{\lambda}^{(\beta+1)}}
\la J_{\lambda+\delta}^{(\beta)},\hGJ J_{\lambda}^{(\beta+1)}\raJ{}
= \frac{\tilc_{\lambda}^{(\beta+1)}}{c_{\lambda}^{(\beta+1)}}
\la J_{\lambda+\delta}^{(\beta)},J_{\lambda+\delta}^{(\beta)}\raJ{}.
\end{eqnarray*}
Applying this relation repeatedly, we have
\[
\la J_{\lambda}^{(\beta)},J_{\lambda}^{(\beta)}\raJ{}
= \prod_{k=0}^{\beta-1}
\frac{\tilc_{\lambda+k\delta}^{(\beta-k)}}{c_{\lambda+k\delta}^{(\beta-k)}}
\la J_{\lambda+\beta\delta}^{(\beta=0)},J_{\lambda+\beta\delta}^{(\beta=0)}
\rangle_{\su}^{(\beta=0)},
\]
which gives the desired result.
\hfill $\Box$\\[2mm]

The norm formula (\ref{eq:NormJack}) can be rewritten into the following
form:
\beq
\la J_{\lambda}^{(\beta)},J_{\lambda}^{(\beta)}\raJ{}
= \frac{(N\beta)!}{(\beta !)^N}
\prod_{(i,j)\in\lambda}
\frac{j-1+\beta(N-i+1)}{j+\beta(N-i)}\cdot
\frac{\lambda_i -j +1 + \beta (\lambda'_j -i)}{
\lambda_i -j + \beta(\lambda'_j -i+1)}.
\label{eq:NormJack2}
\eeq
A proof of the equivalence between (\ref{eq:NormJack}) and 
(\ref{eq:NormJack2}) is given in Appendix B.

\subsection{Shift operators for the multivariable Hermite polynomials}
\label{sec:shiftH}
In this section, we will construct shift operators for the 
multivariable Hermite polynomials.
It should be noted that Heckman has constructed shift operators 
for the Hamiltonian $H_A$ without harmonic potential \cite{H1}.
However, for the application to norm formulas, it is needed to 
compute actions of the shift shift operators on polynomials 
explicitly. Our method gives an unified and straightforward way 
to compute such actions.

To construct shift operators, we first introduce $\syH$ and $\hsyH$ as
follows:
\begin{eqnarray*}
\syH &=& \rhoA(\syJ) = \prod_{i<j}(\beta -\tilhA_i + \tilhA_j),\\
\hsyH &=& \rhoA(\hsyJ) = \prod_{i<j}(-\beta -\tilhA_i + \tilhA_j).
\end{eqnarray*}
Using these operators, we define shift operators for the multivariable
Hermite polynomials:
\[
\GH = \sx^{-1}\syH, \qquad \hGH = \hsyH\sx.
\]
We stress that we have used same $\sx$ as (\ref{eq:sxsy}), and 
therefore $\GH\neq\rhoA(\GJ)$, $\hGH\neq\rhoA(\hGJ)$. 
This reflects the characteristics of the scalar products 
(\ref{eq:spJ}), (\ref{eq:spH}).

If we apply $\rhoA$ to (\ref{eq:sDD}), we have
\begin{eqnarray*}
(s_j+1)(-\beta -\tilhA_j+\tilhA_{j+1}) &=& 
(-\beta -\tilhA_{j+1}+\tilhA_j)(s_j-1),\nonumber\\
(s_j-1)(\beta -\tilhA_j+\tilhA_{j+1}) &=& 
(\beta -\tilhA_{j+1}+\tilhA_j)(s_j+1),\label{eq:shh}\\
s_j (c -\tilhA_j+\tilhA_k)(c -\tilhA_{j+1}+\tilhA_k)
&=& (c -\tilhA_j+\tilhA_k)(c -\tilhA_{j+1}+\tilhA_k) s_j, \nonumber
\end{eqnarray*}
with $c$ arbitrary constant and $k\neq j,j+1$.
These relations imply that 
\beq
\label{eq:yyH}
\syH\left(\cx^{\sSn}\right) = \cx^{-\sSn}, \qquad
\hsyH\left(\cx^{-\sSn}\right) = \cx^{\sSn}.
\eeq
Furthermore, if we apply $\sigA$ to (\ref{eq:actJ}), we see that
\beq
\label{eq:actH}
\syH H_{\lambda+\delta}^{(\beta)} = 
c_{\lambda}^{(\beta+1)}x^{\lambda+\delta}_{w_0} 
    +\mbox{ ``lower terms'' with respect to $\prec$}.
\eeq

For the proof of the shift relations for the Jack polynomials
(Proposition \ref{prop:SrelJ}),
Theorem \ref{thm:dualityJ} played a crucial role.
Here we state analogous result for Hermite case:
\begin{thm}
\label{thm:dualityH}
For $f, g \in \cx^{\sSn}$, 
$\displaystyle
\la \GH f, g\raH{+1} = \la f, \hGH g\raH{}$.
\end{thm}
{\it Proof.}
The proof is similar to that of Theorem \ref{thm:dualityJ}.
For $f,g \in \cx^{\sSn}$, we know that
\[
\la \GH f, g \raH{+1} = 
\int_{-\infty}^{\infty}\cdots\int_{-\infty}^{\infty}
\cP_-(\syH f) g \sx(\GSa)^2\dd x_1\cdots\dd x_N.
\]
On the other hand, applying $\rhoA$ to Lemma \ref{lem:PY2}, we obtain
\[
\cP_- (\syH - \hsyH) = \sum_j \hat{g}_j
(\tilAdag_1,\ldots,\tilAdag_N)(s_j-1)
\quad\mbox{for some}\;\; \hat{g}_j(x_1,\cdots,x_N)\in\cx.
\]
From this relation, we find that 
$\cP_- (\syH - \hsyH) f =0$ for all $f\in\cx^{\sSn}$.
Hence we can replace $\syH$ by $\hsyH$:
\begin{eqnarray*}
\la \GH f, g \raH{+1} &=& 
\int_{-\infty}^{\infty}\cdots\int_{-\infty}^{\infty}
\cP_-(\hsyH f) g \sx(\GSa)^2 \dd x_1\cdots\dd x_N\\
&=& \la \hsyH f, \sx g \raH{} = \la f, \hGH g \raH{}
\end{eqnarray*}
In the last equality, we have used the self-adjointness of the
operator $\tilhA_j$.
\hfill $\Box$\\[2mm]

Now we are in position to state that:
\begin{prop}
\label{prop:SrelH}
$\GH H_{\lambda+\delta}^{(\beta)} = 
c_{\lambda}^{(\beta+1)} H_{\lambda}^{(\beta+1)}$,
$\hGH H_{\lambda}^{(\beta+1)} = \tilc_{\lambda}^{(\beta+1)} 
H_{\lambda+\delta}^{(\beta)}$,
where the constants $c_{\lambda}^{(\beta)}$ and 
$\tilc_{\lambda}^{(\beta)}$ are defined in Lemma \ref{lem:actGJ}.
\end{prop}
{\it Proof.} 
From (\ref{eq:yyH}) and (\ref{eq:actH}), we know that 
$(c_{\lambda}^{(\beta+1)})^{-1}\GH H_{\lambda+\delta}^{(\beta)}$ 
satisfies the first condition of Definition \ref{def:multiH}. 
So it suffice to prove the orthogonality
which can be shown in the same way as Proposition \ref{prop:SrelJ}.
The second equation can be proved in similar way.
\hfill $\Box$\\[2mm]

Using Proposition \ref{prop:SrelH} and Theorem \ref{thm:dualityH},
we can prove the norm formula for $H_{\lambda}^{(\beta)}$.
\begin{prop}[\cite{BF1,vD}]
\label{prop:normH}
\beqa
\lefteqn{
\la H_{\lambda}^{(\beta)}, H_{\lambda}^{(\beta)}\raH{}
= \frac{\pi^{N/2}N!}{2^{|\lambda|+\beta{N(N-1)/2}}}
}\nonumber\\
&& \times \prod_{j=1}^N (\lambda_j + \beta(N-j))! 
\prod_{k=1}^{\beta}\prod_{i<j}
\frac{\lambda_i - \lambda_j -k + \beta(j-i+1)}{
\lambda_i - \lambda_j +k + \beta(j-i-1)},
\label{eq:NormH}
\eeqa
where $|\lambda|=\sum_j\lambda_j$.
\end{prop}
{\it Proof.} 
Using Proposition \ref{prop:SrelH} and Theorem \ref{thm:dualityH},
we see that
\[
\la H_{\lambda}^{(\beta+1)},H_{\lambda}^{(\beta+1)}\raH{+1}
= \frac{\tilc_{\lambda}^{(\beta+1)}}{c_{\lambda}^{(\beta+1)}}
\la H_{\lambda+\delta}^{(\beta)},H_{\lambda+\delta}^{(\beta)}\raH{}.
\]
On the other hand, 
since $H_{\lambda}^{(\beta=0)}(x)$ is direct product of the 
(one-variable) Hermite polynomials, one can evaluate the norm easily:
\[
\la H_{\lambda}^{(\beta=0)},H_{\lambda}^{(\beta=0)}\raH{=0}
= \frac{\pi^{N/2}\cdot\#\Sn^{\lambda}}{
2^{|\lambda|}}\prod_{j=1}^N \lambda_j ! .
\]
Using these relations, one arrives at the formula above.
\hfill $\Box$\\[2mm]

The norm formula (\ref{eq:NormH}) can be rewritten into the following
form \cite{BF1}:
\begin{eqnarray*}
\lefteqn{
\la H_{\lambda}^{(\beta)},H_{\lambda}^{(\beta)}\raH{}
= \frac{\pi^{N/2}}{2^{|\lambda|+\beta N(N-1)/2}}
\cdot\frac{\prod_{j=1}^N (j\beta)!}{(\beta !)^N}
}\nonumber\\
&&\times\prod_{(i,j)\in\lambda}
\frac{\{j-1+\beta(N-i+1)\}
      \{\lambda_i -j +1 + \beta (\lambda'_j -i)\}}{
      \lambda_i -j + \beta(\lambda'_j -i+1)}.
\label{eq:NormH2}
\end{eqnarray*}
It should be remarked that other proofs of these formulas have been 
given via limiting procedure \cite{BF1,vD}.

\subsection{Shift operators for the multivariable Laguerre polynomials}
We first define $\sxL$, $\syL$ and $\hsyL$ as follows:
\begin{eqnarray*}
\sxL &=& \prod_{i<j}(z_i^2-z_j^2),\\
\syL &=& \rhoB(\syJ) = \prod_{i<j}(\beta -\tilhB_i/2 + \tilhB_j/2),\\
\hsyL &=& \rhoB(\hsyJ) = \prod_{i<j}(-\beta -\tilhB_i/2 + \tilhB_j/2).
\end{eqnarray*}
After same discussion as in the previous subsection, we see that
\beq
\label{eq:yyL}
\syL\left(\czz^{\sSn}\right) = \czz^{-\sSn}, \qquad
\hsyL\left(\czz^{-\sSn}\right) = \czz^{\sSn},
\eeq
and
\beq
\label{eq:actL}
\syH L_{\lambda+\delta}^{(\beta)} = 
c_{\lambda}^{(\beta+1)}z^{2(\lambda+\delta)}_{w_0} 
    +\mbox{ ``lower terms'' with respect to $\prec$}.
\eeq

Now we define the shift operators for Laguerre case:
\[
\GL = \sxL^{-1}\syL, \qquad \hGL = \hsyL\sxL.
\]
These operators enjoy the following properties:
\begin{thm}
\label{thm:dualityL}
For $f, g \in \czz^{\sSn}$, 
$\displaystyle
\la \GL f, g\raL{+1} = \la f, \hGL g\raL{}$.
\end{thm}
{\it Proof.}
The proof for this case is also similar to that of 
Theorem \ref{thm:dualityJ}.
For $f,g \in \czz^{\sSn}$, we know that
\[
\la \GL f, g \raL{+1} = 
\int_{-\infty}^{\infty}\cdots\int_{-\infty}^{\infty}
\cP_-(\syL f) g \sxL(\GSb)^2 \dd z_1\cdots\dd z_N.
\]
On the other hand, applying $\rhoB$ to Lemma \ref{lem:PY2}, we find that
\[
\cP_- (\syL - \hsyL) = \sum_j \hat{g}_j
(\tilBdag_1,\ldots,\tilBdag_N)(s_j-1)
\quad\mbox{for some}\;\; \hat{g}_j(x_1,\cdots,x_N)\in\cx.
\]
From this relation, we see that 
$\cP_- (\syL - \hsyL) f =0$ for all $f\in\czz^{\sSn}$.
Hence we can replace $\syL$ by $\hsyL$:
\begin{eqnarray*}
\la \GH f, g \raL{+1} &=& 
\int_{-\infty}^{\infty}\cdots\int_{-\infty}^{\infty}
\cP_-(\hsyL f) g \sxL(\GSa)^2 \dd z_1\cdots\dd z_N\\
&=& \la \hsyL f, \sxL g \raL{} = \la f, \hGL g \raL{}
\end{eqnarray*}
In the last equality, we have used the self-adjointness of 
the operator $\tilhB_j$.
\hfill $\Box$\\[2mm]

We then state the following results:
\begin{prop}
\label{prop:SrelL}
$\GL L_{\lambda+\delta}^{(\beta)} = 
c_{\lambda}^{(\beta+1)} L_{\lambda}^{(\beta+1)}$,
$\hGL L_{\lambda}^{(\beta+1)} = \tilc_{\lambda}^{(\beta+1)} 
L_{\lambda+\delta}^{(\beta)}$,
where the constants $c_{\lambda}^{(\beta)}$ and 
$\tilc_{\lambda}^{(\beta)}$ are defined in Lemma \ref{lem:actGJ}.
\end{prop}
{\it Proof.} 
From (\ref{eq:yyL}) and (\ref{eq:actL}), we know that 
$(c_{\lambda}^{(\beta+1)})^{-1}\GL L_{\lambda+\delta}^{(\beta)}$ 
satisfies the first condition of 
Definition \ref{def:multiL} up to a constant factor. 
So it suffice to prove the orthogonality
which can be shown in the same way as Proposition \ref{prop:SrelJ}.
The second equation can be proved in similar way.
\hfill $\Box$\\[2mm]

Using Proposition \ref{prop:SrelL} and Theorem \ref{thm:dualityL},
we can prove the norm formula for $L_{\lambda}^{(\beta)}$.
\begin{prop}[\cite{BF1,vD}]
\beqa
\lefteqn{
\la L_{\lambda}^{(\beta)}, L_{\lambda}^{(\beta)}\raL{}
= N! \prod_{j=1}^N (\lambda_j + \beta(N-j))! 
}\nonumber\\
&&\times
\prod_{j=1}^N \Gamma(\lambda_j +\beta(N-j) + \gamma + 1/2)
\prod_{k=1}^{\beta}\prod_{i<j}
\frac{\lambda_i - \lambda_j -k + \beta(j-i+1)}{
\lambda_i - \lambda_j +k + \beta(j-i-1)}
\label{eq:NormL}
\eeqa
where $\Gamma(\cdot)$ denotes the gamma function.
\end{prop}
{\it Proof.} 
The proof of this proposition is similar to the Hermite case.
We only note the following formula for the case $\beta=0$:
\[
\la L_{\lambda}^{(\beta=0)},L_{\lambda}^{(\beta=0)}\raH{=0}
= (\#\Sn^{\lambda})\prod_{j=1}^N \left\{\lambda_j ! 
\cdot\Gamma(\lambda_j+\gamma+1/2)\right\} ,
\]
which follows from the norm formula of the one-variable 
Laguerre polynomials.
\hfill $\Box$\\[2mm]

The norm formula (\ref{eq:NormL}) can be rewritten into the following
form \cite{BF1}:
\begin{eqnarray*}
\lefteqn{
\la L_{\lambda}^{(\beta)},L_{\lambda}^{(\beta)}\raH{}
= \frac{\prod_{j=1}^N (j\beta)!}{(\beta !)^N}
\prod_{j=1}^N \Gamma(\lambda_j +\beta(N-j) + \gamma + 1/2)
}\nonumber\\
&&\times\prod_{(i,j)\in\lambda}
\frac{\{j-1+\beta(N-i+1)\}
      \{\lambda_i -j +1 + \beta (\lambda'_j -i)\}}{
      \lambda_i -j + \beta(\lambda'_j -i+1)}.
\label{eq:NormL2}
\end{eqnarray*}
It should be remarked that other proofs of these formulas have been 
given via limiting procedure \cite{BF1,vD}.

\section{Concluding remarks}
In this paper, we have constructed the intertwining operators 
that map the Jack polynomials to the multivariable Hermite and 
Laguerre polynomials.

We restrict ourselves to symmetric polynomials though the operators
$\sigA$ and $\sigB$ are applicable to non-symmetric case, i.e.
we can obtain the non-symmetric counterparts of the multivariable 
Hermite and Laguerre polynomials:
\[
E^{(\mbox{\rm\scriptsize H})\lambda}_w (x)
= 2^{-|\lambda|/2}E^{\lambda}_w (\tilAdag)\cdot 1, \qquad
E^{(\mbox{\rm\scriptsize L})\lambda}_w (z)
= E^{\lambda}_w ((\tilBdag)^2/2)\cdot 1.
\]
Baker and Forrester named these polynomials non-symmetric Hermite 
and Laguerre polynomials respectively, and studied their properties
\cite{BF2,BF3}.
We note that some of their results may be obtained directly form
the corresponding properties of the Jack polynomials by applying
the intertwiners.

Our constructs are based on the degenerate double affine Hecke 
algebra, so it is expected that the results given here extend to 
non-degenerate, i.e. $q$-deformed case.
As van Diejen \cite{vD} already proposed $q$-difference counterpart of the 
Hamiltonians $\Ha$ and $\Hb$, it would be nice to
clarify algebraic structure of the $q$-cases.
We hope to report on them in the near future.

\section*{Appendices}
\addcontentsline{toc}{chapter}{\protect\numberline{}{Appendices}}
\subsection*{Appendix A: Proof of Lemma 4.3}
\renewcommand{\theequation}{\mbox{A.}\arabic{equation}}
\setcounter{equation}{0}
\renewcommand{\thelem}{\mbox{A.}\arabic{lem}}
\setcounter{lem}{0}
In Appendix A, we will give a proof of Lemma \ref{lem:PY2}.
We remark again that the proof given in this section is
a limiting case of \cite{KJr}.

We begin with seeing some properties of the anti-symmetrizer.
\begin{lem}[\cite{KJr}]
\label{lemma:divisibility}
\begin{enumerate}
\item
The anti-symmetrizer $\cP_-$ is divisible by $1+(-1)^{l(w_0)}w_0$
both on the left and on the right.
\item 
For all $j=1,\ldots,N-1$,
the anti-symmetrizer $\cP_-$ is divisible by $1-s_j$ 
both on the left and on the right.
\end{enumerate}
\end{lem}
{\it Proof.} 
(i) $\Sn$ can be divided into pairs $(w,w w_0)$. 
Then, rewriting into the summation over such pairs, we have
\begin{eqnarray*}
\cP_- &=& 
\sum_{(w,w w_0)}\left\{(-1)^{l(w)}w + (-1)^{l(ww_0)}ww_0\right\}\\
&=& \sum_{(w,w w_0)}(-1)^{l(w)}w\left\{1+(-1)^{l(w_0)}w_0\right\}.
\end{eqnarray*}
Divisibility on the left is proved similarly.

(ii) Can also be proved by similar discussion.
\hfill $\Box$\\[2mm]

From Lemma \ref{lemma:divisibility} (ii), we know that 
$\Ker\cP_- \supset \sum_j \Ker (1-s_j)$. 
To describe kernel of the anti-symmetrizer, we first investigate 
kernels of $1-s_j$ and their union. 
\begin{lem}[\cite{KJr}]
\label{lemma:Ker}
\begin{enumerate}
\item 
Let $V$ is a representation of $\Sn$, and denote 
$V_j = \Ker (1-s_j)$, $V'=\sum_j V_j$. Then 
$V'$ is $\Sn$-invariant.
\item 
Assume $V$ is a finite-dimensional irreducible representation
of $\Sn$. Then we have
\[
V'= \left\{\begin{array}{ll}
0\quad  \mbox{(if } V \mbox{ is the sign representation),}\\
V\quad \mbox{(otherwise)}.
\end{array}\right.
\]
\end{enumerate}
\end{lem}
{\it Proof.} 
(i) 
From the definition of $V_j$, it follows that $s_j (s_i v)=s_i v$
for all $v\in s_i V_j$. 
If we introduce $v_{\pm}=(v\pm s_iv)/2$,
we see that $s_j (v_+ - v_-)= v_+ - v_-$ 
which means $v_+ - v_-\in V_j$. 
Since $v_+\in V_i$ by definition, we obtain
$v=v_+ + v_- \in V_i + V_j$.
This leads to $s_i V_j\subset V_i+V_j$, which concludes the proof.

(ii) 
From (i), it follows that $V'$ is a subrepresentation.
Due to the irreducibility, $V'$ can be either $0$ or $V$.
If $V'=0$, then we have $V_j =0$ for all $j$. 
This means that $1-s_j$ is invertible,
i.e. for all $v\in V$, there exists $u$ such that $v=(1-s_j)u$.
Then we obtain $s_j v=-v$ for all $v\in V$, i.e. $V$ is the sign
representation.\hfill $\Box$\\[2mm]%
From Lemma \ref{lemma:Ker}(ii), it immediately follows:
\[
\Ker\cP_- = \sum_j \Ker (1-s_j)
\]
for any finite-dimensional representation of $\Sn$. 
Note this identity also holds for the representation 
of $\Sn$ in the space of polynomials $\cx$, since this representation 
is a direct sum of finite-dimensional representations. 

We now introduce operators $\hs_j$ as
\[
\hs_j = s_j + \beta\frac{s_j -1}{x_j - x_{j+1}}.
\]
Using these operators, we can define another representation of 
the degenerate affine Hecke algebra $\daH$ on $\cx$:
\[
\rho' (\hDA_j) = x_j, \qquad 
\rho' (s_j) = \hs_j.
\]
Using the isomorphism $\rho'$, we introduce deformed 
anti-symmetrizer $\cPb$ as 
\[
\cPb = \rho'(\cP_-) = \frac{1}{\# \Sn}\sum_{w\in \sSn}(-1)^{l(w)}\hw,
\]
where $\hw = \rho' (w)$.
\begin{lem}[\cite{KJr}]
\label{lem:cPb}
$\Ker \cPb = \Ker \cP$ for the action of $\cPb$ in $\cx$:
\end{lem}
{\it Proof.} 
By similar discussion to Proposition \ref{lemma:divisibility},
we know that
$\cPb$ is divisible by $\hs_j-1$ both on the left and on the right
for every $j=1,\ldots,N-1$.
Hence we have
\beq
\label{eq:supset}
\Ker \cPb \supset \sum_j \Ker (1-\hs_j) 
= \sum_j \Ker (1-s_j) = \Ker \cP,
\eeq
and thus $\dim(\Ker\cPb)\geq\dim(\Ker\cP_-)$.
On the other hand, if we denote $\cx_n$ as space of polynomials of
order $n$, it is clear that $\cPb$ preserves $\cx_n$.
Since $\dim(\Ker\cPb)$ can not decrease under specialization,
it follows that 
$\dim(\Ker\cPb) \leq \dim(\Ker\cP_-^{(\beta=0)})=\dim(\Ker\cP_-)$ 
and hence we have 
\beq
\label{eq:dimKer}
\dim(\Ker\cPb) = \dim(\Ker\cP_-).
\eeq
Thus it follows from (\ref{eq:supset}) and (\ref{eq:dimKer}) 
that $\Ker\cPb=\Ker\cP=\sum_j\Ker(1-s_j)$.
\hfill $\Box$\\[2mm]

We then define $\sy'$ and $\hsy'$ as
\begin{eqnarray*}
\sy' &=& \rho'(\syJ) = \prod_{i<j}(\beta -x_i + x_j),\\
\hsy' &=& \rho'(\hsyJ) = \prod_{i<j}(-\beta -x_i + x_j).
\end{eqnarray*}
\begin{lem}[\cite{KJr}]
\label{lem:cpY}
\[
\cP_- (\sy' - \hsy') f = 0 \quad\mbox{for all}\quad f\in\cx^{\sSn}.
\]
\end{lem}
{\it Proof.} 
We can show that
\[
(1+(-1)^{N(N-1)/2}w_0)(\sy'-\hsy') 
= (\sy'-\hsy')(1-w_0).
\]
by the direct calculation. 
Considering the action on $\cx^{\sSn}$, we have
\[
(1+(-1)^{N(N-1)/2}w_0)(\sy'-\hsy') f =0,
\]
for all $f\in\cx^{\sSn}$.
Using this formula and Proposition \ref{lemma:divisibility} (i),
we obtain the desirous result.
\hfill $\Box$\\[2mm]

From Lemma \ref{lem:cPb} and Lemma \ref{lem:cpY}, we know that
\[
\label{eq:cPbY}
\cPb (\sy' - \hsy') f = 0 \quad\mbox{for all}\;\; f\in\cx^{\sSn}.
\]
On the other hand, the following statement can easily be proved:
\begin{lem}
\label{lemma:DivThm}
Let $\hat{A}$ be a operator of the form,
 $\hat{A}=\sum_{w\in \sSn}g_w\hw$ with $g_w\in\cx$.
If $\hat{A}f =0$ for all $f\in\cx^{\sSn}$,
then $\hat{A}$ can be represented in the following form:
\[
\hat{A} = \sum_j \hat{g}_j (\hs_j-1)
\quad\mbox{for some}\;\; \hat{g}_j\in\cx^{\sSn}.
\]
\end{lem}
{\it Proof.} 
The operator $\hat{A}$ can be rewritten as
\[
\hat{A} = 
\sum_{j_1,\ldots,j_k}\hat{g}'_{j_1,\ldots,j_k}
(\hs_{j_1}-1)\cdots(\hs_{j_k}-1) + \hat{g}'_0
\quad\mbox{for some}\;\; \hat{g}'_{j_1,\ldots,j_k}\in\cx.
\]
Then the assumption of the proposition means $\hat{g}'_0=0$,
which gives the desirous result.
\hfill $\Box$\\[2mm]

Applying Lemma \ref{lemma:DivThm} to Lemma \ref{lem:cpY}, 
we conclude that 
\beq
\label{eq:PY1}
\cPb (\sy' - \hsy') = \sum_j \hat{g}_j(x_1,\cdots,x_N)(\hs_j-1)
\quad\mbox{for some}\;\; \hat{g}_j(x_1,\ldots,x_N)\in\cx.
\eeq
Applying $(\rho')^{-1}$ completes the proof of Lemma 4.3.

\subsection*{Appendix B: Equivalence of the two expressions for the 
norm formula}
\renewcommand{\theequation}{\mbox{B-}\arabic{equation}}
\setcounter{equation}{0}
\renewcommand{\thelem}{\mbox{B}\arabic{lem}}
\setcounter{lem}{0}
In Appendix B, we will give a proof of equivalence between two
expressions of the norm formulas.
We first begin with considering the Jack case.

Let $\lambda$ be a partitions satisfying the following conditions
(see Figure 1 below):
\begin{eqnarray*}
\lefteqn{\lambda_{p-1} > \lambda_p = \cdots = \lambda_{p+r_1-1}
> \lambda_{p+r_1} = \cdots = \lambda_{p+r_1+r_2-2}}\\
&& > \cdots > 
\lambda_{p+r_1+\cdots+r_{m-1}} = \cdots = \lambda_{p+r_1+\cdots+r_m-1}
> \lambda_{p+r_1+\cdots+r_m} = \cdots =0,\\
\lefteqn{\lambda'_1 = \cdots = \lambda'_{s_1}
> \lambda'_{s_1+1} = \cdots = \lambda'_{s_1+s_2}}\\
&& > \cdots >
\lambda'_{s_1+\cdots+s_{m-1}+1} = \cdots = \lambda'_{s_1+\cdots+s_m}
> \lambda'_{s_1+\cdots+s_m} = \cdots = 0.
\end{eqnarray*}
We further define $\mu$ as
$\mu = (\lambda_1,\cdots,\lambda_p+1,\cdots,\lambda_N)$.

\begin{center}
{\leavevmode\epsfbox{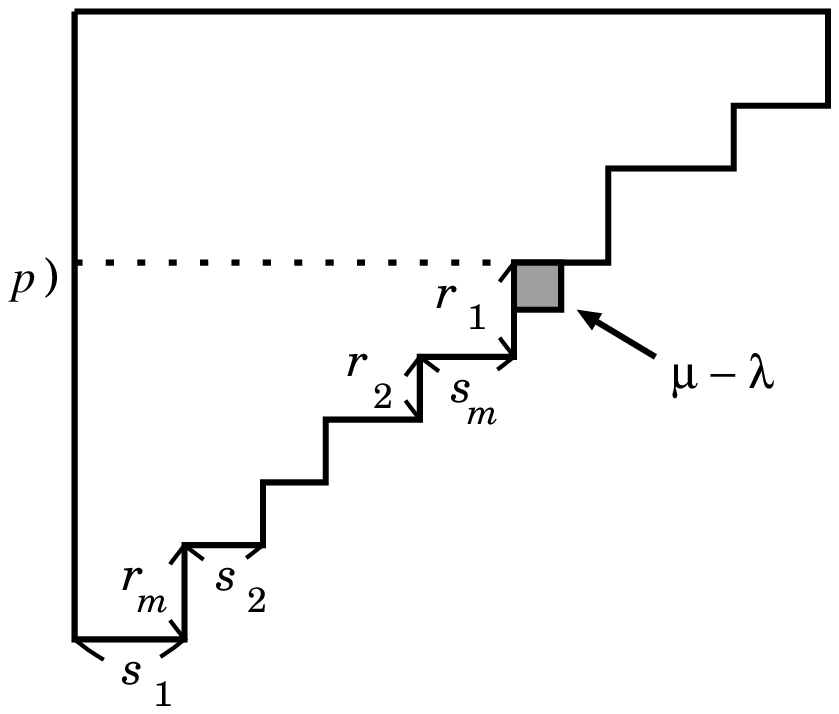}}
{\bf Figure 1:} Young diagram of $\lambda$ and $\mu$
\end{center}

Calculating the ratio 
$\la J^{(\beta)}_{\mu}, J^{(\beta)}_{\mu}\raJ{}/
\la J^{(\beta)}_{\lambda}, J^{(\beta)}_{\lambda}\raJ{}$
by using (\ref{eq:NormJack}) or (\ref{eq:NormJack2}),
one can show that both cases reduce to 
\begin{eqnarray*}
\lefteqn{
\frac{\la J^{(\beta)}_{\mu}, J^{(\beta)}_{\mu}\raJ{}}{
\la J^{(\beta)}_{\lambda}, J^{(\beta)}_{\lambda}\raJ{}}
= \prod_{i=1}^{p-1}
\frac{\lambda_i-\lambda_p+\beta(p-i)}{\lambda_i-\lambda_p+\beta(p-i-1)}
\cdot\frac{\lambda_i-\lambda_p-1+\beta(p-i)}{
\lambda_i-\lambda_p-1+\beta(p-i-1)}}\\
&&\times \frac{s_m+1+\beta(r_1-1)}{1+\beta(r_1-1)}
\cdots
\frac{s_m+\cdots+s_1+1+\beta(r_1+\cdots+r_m-1)}{
s_m+\cdots+s_2+1+\beta(r_1+\cdots+r_m-1)}\\
&&\times \frac{\beta r_1}{s_m+\beta r_1}\cdot
\frac{s_m+\beta(r_1+r_2)}{s_m+s_{m-1}+\beta(r_1+r_2)}
\cdots
\frac{s_m+\cdots s_2+\beta(r_1+\cdots+r_m)}{
s_m+\cdots+s_1+\beta(r_1+\cdots+r_m)}\\
&&\times \frac{s_m+\cdots s_1+\beta(N-p+1)}{
s_m+\cdots+s_1+1+\beta(N-p)}
\cdot\frac{1}{\beta}.
\end{eqnarray*}

On the other hand, if we consider the simplest case $\lambda=\phi$,
both (\ref{eq:NormJack}) and (\ref{eq:NormJack2})
reduce to $\la 1,1 \raJ{} = (\beta N)!/(\beta !)^N$.
Hence, by induction, we conclude that (\ref{eq:NormJack}) and 
(\ref{eq:NormJack2}) are equivalent for all $\lambda$.

The Hermite and Laguerre cases can be proved in the similar fashion.

\end{document}